\documentclass[11 pt, psamsfonts]{amsart}
\usepackage{amssymb,color}

\def\ignore#1{\relax}

 \definecolor{red}{rgb}{1,0,0}
 \definecolor{blue}{rgb}{.1,0,1}
 \DeclareMathOperator{\FPdim}{FPdim}
  
\def\g{\mathfrak g}

\def\h{\mathfrak h}

\def\sl{\mathfrak{sl}}
\def\so{\mathfrak{so}}
\def\sp{\mathfrak{sp}}
\def\Vl{{M_\la}}
\def\R{{\mathbb R}}
\def\Z{{\mathbb Z}}

\def\Q{{\mathbb Q}}
\def\C{{\mathbb C}}

\def\la{\lambda}
\def\ep{\epsilon}

\def\N{\mathbb N}

\def\A{\mathcal A}
\def\Ca{\mathcal C}

\def\one{\mathbf 1}

\def\A{\mathcal A}

\def\Da{\mathcal D}

\def\i{{\rm i}}

\def\B{{\mathcal B}}
\def\U{{\bf U}}

\def\ignore#1{\relax}
\def\om{\omega}
\def\e{\epsilon}
\def\1{{\bf 1}}

\def\End{{\rm End}}
\def\Hom{{\rm Hom}}
\def\ve{\varepsilon}

\def\ot{\otimes}
\def\Spec{{\rm Spec}}
\def\Rep{{\rm Rep}}

\calclayout

\renewenvironment{proof}[1][\proofname]{\par
  \normalfont
  \topsep6\p@\@plus6\p@ \trivlist
  \item[\hskip\labelsep\bfseries
    #1\@addpunct{.}]\ignorespaces
}{%
  \qed\endtrivlist
}

\def\th@plain{%
  \let\thmhead\thmhead@plain \let\swappedhead\swappedhead@plain
  \thm@preskip.5\baselineskip\@plus.2\baselineskip
                                    \@minus.2\baselineskip
  \thm@postskip\thm@preskip
  \itshape
\renewcommand{\labelenumi}{{(\alph{enumi})\quad}}
                        \renewcommand{\labelenumii}{{(\roman{enumii})\ }}
}
\def\th@definition{%
  \let\thmhead\thmhead@plain \let\swappedhead\swappedhead@plain
  \thm@preskip.5\baselineskip\@plus.2\baselineskip
                                    \@minus.2\baselineskip
  \thm@postskip\thm@preskip
  \upshape
}
\def\th@remark{%
  \thm@headfont{\itshape}
  \let\thmhead\thmhead@plain \let\swappedhead\swappedhead@plain
  \thm@preskip.5\baselineskip\@plus.2\baselineskip
                                    \@minus.2\baselineskip
  \thm@postskip\thm@preskip
  \upshape
}

{\theoremstyle{plain}
\newtheorem{theorem}{Theorem}[section]
}

{\theoremstyle{plain}
\newtheorem{proposition}[theorem]{Proposition}
}

{\theoremstyle{plain}
\newtheorem{corollary}[theorem]{Corollary}
}

{\theoremstyle{plain}
\newtheorem{lemma}[theorem]{Lemma}
}

{\theoremstyle{plain}
\newtheorem{conjecture}[theorem]{Conjecture}
}

{\theoremstyle{definition}

}

{\theoremstyle{definition}

}

{\theoremstyle{remark}
\newtheorem{remark}[theorem]{Remark}
}

{\theoremstyle{remark}

}

\numberwithin{equation}{section}

\renewcommand{\labelenumi}{{ \theenumi.}}
\renewcommand{\labelenumii}{{(\alph{enumii})}}

\def\la{\lambda}
\def\al{\alpha}

\def\choose #1 #2{\begin{pmatrix}#1\\#2\end{pmatrix}}

\def\i{{\rm i}}

\def\Aut{{\rm Aut}}

\begin{document}

\title{$SO(N)_2$ Braid group representations are Gaussian}

\author{Eric C. Rowell}
\address{Department of Mathematics\\ Texas A\&M University\\ College Station,
Texas}
\email{rowell@math.tamu.edu}

\author{Hans Wenzl}

\address{Department of Mathematics\\ University of California\\ San Diego,
California}

\email{hwenzl@ucsd.edu}

\thanks{E.C.R. is partially supported by US NSF grant 1108725, and wishes to thank V.F.R. Jones, Z. Wang, D. Naidu and C. Galindo
for illuminating discussions.  Part of this was written while E.C.R. was visiting BICMR, Peking University and the Chern Institute, Nankai University
and he gratefully acknowledges the hospitality of these Institutions.}

\begin{abstract} We give a description of the centralizer algebras for tensor powers of spin objects in the pre-modular categories
$SO(N)_2$ (for $N$ odd) and $O(N)_2$ (for $N$ even) in terms of quantum $(n-1)$-tori, via non-standard deformations of $U\so_N$.
As a consequence we show that the corresponding braid group representations are Gaussian representations,
the images of which are finite groups.
This verifies special cases of a conjecture that braid group
representations coming from weakly integral braided fusion
categories have finite image.
\end{abstract}\maketitle

\section{Introduction}
Let $V$ be the $N$-dimensional vector
representation of the quantum group $U_q\g$, where $\g\in\{\sl_N,\so_N,\sp_N\}$
and where $q$ is an indeterminate. Let $\B_n$ be the braid group on $n$-strands, for a natural number $n>0$. Then the centralizer algebras $\End_{U_q\g}(V^{\ot n})$ each have complete
descriptions in terms of semisimple quotients of braid group algebras $\C(q)\B_n$, namely Hecke and
BMW-algebras (\cite{Jimbo}, \cite{bcd}).  For $q=e^{\pi \i/\ell}$ the representation categories $Rep(U_q \g)$ are not semisimple, but have semisimple sub-quotients $\Ca(\g,\ell)$ via a process called ``purification'' in \cite{Turaev}. Continuing to denote by $V$ the image of $V$ in the sub-quotient $\Ca(\g,\ell)$, the centralizer algebras $\End({V}^{\ot n})$ in the fusion category $\Ca(\g,\ell)$ still are quotients of Hecke or BMW-algebras, so that the description in terms of the braid group algebra persists.  The (closure) of the image of these braid group representations were analyzed in \cite{FLW},\cite{RLW}, which provided evidence for the following conjecture (see \cite{NR,RW}):
\begin{conjecture}\label{mainconj}
Let $\Ca$ be a braided fusion category and let $X$ be a simple object in $\Ca$.  The braid group representations $\B_n$ on $\End(X^{\ot n})$
have finite image for all $n>0$ if and only if $\FPdim(X)^2\in\Z$.
\end{conjecture}
Here $\FPdim(X)\in\R$ is the Frobenius-Perron dimension, which coincides with the categorical dimension for unitary fusion categories.
Categories with $\FPdim(X)^2\in\Z$ for all simple objects $X$ are called \emph{weakly integral}.  One large class of weakly integral braided fusion categories for which this conjecture has been verified are the so-called \emph{group theoretical} categories related to the (doubles of) finite groups (see \cite{ERW}).  The main difficulty in verifying Conjecture \ref{mainconj} for arbitrary objects $X$ in a braided fusion category is that a sufficiently explicit description of the braid group action on $\End(X^{\ot n})$ is usually lacking.

For the class of categories $\Ca(\g,\ell)$ associated with quantum groups
at roots of unity
the ``only if" part of the
conjecture has been confirmed (see \cite{Rarg}); more precisely, it was shown that for any object $X$ in a fusion category $\Ca(\g,\ell)$ for which $\FPdim(X)^2\not\in\Z$ the associated braid representations have infinite image. Conversely, for simple objects $X$ in $\Ca(\g,\ell)$ for which $\FPdim(X)^2\in\Z$ the only remaining open cases are in $\Ca(\so_N,2N)$ for $N$ odd and $\Ca(\so_N,N)$ for $N$ even.  We will adopt the uniform notation $SO(N)_2$ for these two families (this notation conforms with the physics literature, where the subscript $2$ is the level).
For $N$ odd, the (fundamental) spinor object $S\in SO(N)_2$ has dimension $\sqrt{N}$, whereas for $N$ even the two (fundamental) spinor objects $S_\pm\in SO(N)_2$ have dimension $\sqrt{N/2}$.  In particular, Conjecture \ref{mainconj} predicts that
the $\B_n$ representations associated with these have finite image.

The main result of our paper is to prove Conjecture \ref{mainconj} for $SO(N)_2$, for all positive integers $N$.
It had been verified for $N\leq 8$ using low rank coincidences and when $\sqrt{N}$ or $\sqrt{N/2}$ are integral, by appealing to the results of \cite{ERW}, see \cite{NR}.
For general $N$, it seems to be difficult to give an intrinsic description
of the braid representations as the number of eigenvalues of the image of
a standard generator is increasing with N.
We overcome these difficulties by the following approach: We denote by $S$ the quantum version of the spinor representation of $U_q\so_N$ (for $N$ odd) as well as its image in $SO(N)_2$ and by $S_{\pm}$ the analogous fundamental spinor objects for $N$ even and their images in $SO(N)_2$.  For $N$ even we consider a semidirect (smash) product $U_q\so_N\rtimes \Z_2$ and we will also denote by $S$ the irreducible $U_q\so_N\rtimes \Z_2$-module whose restriction to $U_q\so_N$ is $S_-\oplus S_+$.  For generic parameter $q$, the centralizer algebras
$\End_\U(S^{\ot n})$  are described (\cite[Theorem 4.8]{WSp}) in terms of a non-standard
deformation $U_q^\prime \so_n$ of $U\so_n$, for both $N$ odd and even.  Although $\Rep(\U)$ carries a braiding, the image of $\B_n$ inside $\End_{\U}(S^{\ot n})$ does not generate these algebras.  On the
other hand, for
$q$ a $2N$th root of unity, we show that the algebra $U_q^\prime
\so_n$ admits a homomorphism into the quantum $(n-1)$-torus
$T_q(n)$, which contains an isomorphic copy of $\End(S^{\ot n})$. The key observation now is that this homomorphism identifies the image of the $\B_n$-representations in $\End(S^{\ot n})$ for the braided fusion category $SO(N)_2$ with
so-called Gaussian braid representations (so named because the coefficients are Gaussian functions of the form $Ke^{2\pi \i a^2/\ell}$, defined in Proposition \ref{identifiedBnreps}(a)) which
live in the quantum torus.  These explicitly realized braid representations can be shown to have
finite images, which implies the conjecture for $SO(N)_2$.  So the identification of these different braid representations is achieved using the representation theory of the algebra $U'_q\so_n$.

Here is a more detailed outline of the contents of this article. In Section \ref{section: duality} we review results
about the centralizer algebras $\End(S^{\otimes n})$ where $S$ is a spinor representation of $U_q\so_N$ respectively $U_q\so_N\rtimes \Z_2$,
or the corresponding object in one of the associated fusion categories. Most of these results have already more or less appeared
before in \cite{Has}, \cite{WSp}. In Section \ref{section: rep. theory}, we reprove and extend several results by Klimyk and his coauthors concerning the
representation theory of $U'_q\so_n$. In contrast with \emph{loc. cit.}, we use a Verma module approach which also has
the advantage of proving (crucial, for our paper) uniqueness results at roots of unity, for certain types of modules.
In Section \ref{section: quantum torus} we construct representations of $U'_q\so_n$ into algebras called quantum tori.  The main result of this section is the identification of these representations with those in $\End(S^{\otimes n})$ for fusion categories $SO(N)_2$ ($N$ odd) and $O(N)_2$ ($N$ even). This allows us to describe the corresponding tower of centralizer algebras in terms of the quantum tori using Jones' basic construction.  Finally, we identify the braid group representations corresponding to the object $S$ in $SO(N)_2$ (respectively $O(N)_2$) for $N$ odd (respectively $N$ even) with the Gaussian braid representations  first encountered in the work of Jones and Goldschmidt \cite{GJ},\cite{jonespjm} for $N$ odd. The easy generalization to $N$ even is worked out in \cite{GR}.
From this one easily verifies Conjecture \ref{mainconj}
in our case.

\section{Duality for spinor representations}\label{section: duality}

\subsection{Deformations of $U\so_n$}
The algebra $U'_q\so_n$ is defined (see \cite{GK})
 via generators $B_1,\ldots, B_{n-1}$ satisfying the relations
$B_iB_j=B_jB_i$ for $|i-j|\neq 1$ and the $q$-Serre relations:
\begin{equation}\label{quantumSerre}
B_i^2B_{i\pm 1}-(q+q^{-1})B_iB_{i\pm 1}B_i +B_{i\pm 1}B_i^2= B_{i\pm 1}.
\end{equation}
It is well-known that in the classical limit $q=1$ we obtain a presentation of the universal enveloping algebra $U\so_n$ of the orthogonal Lie algebra $\so_n$,
and for this reason $U'_q\so_n$ is sometimes called the \emph{non-standard} deformation of $U\so_n$.
It follows from the definitions that the elements $B_1, B_3, \ldots, B_{n-1}$ (for $n$ even, respectively $B_{n-2}$ for $n$ odd), generate an abelian subalgebra $A$  of $U_q'\so_n$.
We define a weight vector of a $U_q'\so_n$-module $V$ to be a common eigenvector of
the generators of $A$. We call a weight regular if all the eigenvalues of generators
$B_{2i-1}$ of $A$ are of the form $[r]$ with $r$ an integer or a half integer,
and $[r]=(q^r-q^{-r})/(q-q^{-1})$ the usual $q$-number.

 In the following we denote by $\U$ a semidirect
product of the (standard) Drinfeld-Jimbo quantum group $U_q\so_N$ with $\Z_2$.   For $N$ odd, $\U$ is just
the direct sum of the corresponding $\C$-algebras, while in the $N$ even case,
the nontrivial element $t$ of $\Z_2$ acts via the obvious type $D_{N/2}$ graph automorphism.  This completely determines the
defining relations for $\U$.  It is also easy to check that the map $\Delta(t)=t\otimes t$ extends the bialgebra  structure
of $U_q\so_N$ to $\U$.  Indeed, by \cite[Theorem 2.1]{M} $\U$ (called the \emph{smash product algebra} in \emph{loc. cit.})
is a ribbon Hopf algebra as the action of $t$ preserves the braiding.
 For $N$ odd, it is clear that $\Rep(\U)\cong \Rep(U_q\so_N)\boxtimes \Rep(\Z_2)$  (Deligne tensor product) as ribbon categories.
 Note that by \cite{GM} $\Rep(\U)$ is the $\Z_2$-equivariantization of $\Rep(U_q\so_N)$.
We shall also be interested in the case where $q$ is a root of unity.  In this case we consider the subcategory of
tilting modules in $\Rep(\U)$ which is again a ribbon category (see e.g. \cite{Wcstar} for details).  As such, we may consider the quotient category
by negligible morphisms (see \cite[Section XI.4]{Turaev}) to obtain ribbon fusion categories $SO(N)_r$ and $O(N)_r$, which we describe below.

The algebra $\U$ should be viewed as a quantum version of $Pin(N)$.  Indeed, $\U$ is well-defined in the classical limit $q=1$,
and its finite dimensional simple representations are in 1-1 correspondence with the simple representations of $Pin(N)$.
It is easy to see that we obtain a well-defined quantum version of the spinor $Pin(N)$-module $S$ for $\U$
(where the matrices of the generators $E_i$, $F_i$ and $t$ do not depend on $q$).
As any finite-dimensional simple $Pin(N)$-module does appear in some tensor power of $S$, we can also make it into a $\U$-module.
This deformation also works for roots of unity. It is well-known
and easy to check that if the restriction of a simple $U_q\so_N\rtimes \Z_2$-module to $U_q\so_N$ does not remain simple, it decomposes into the direct sum of two irreducible $U_q\so_N$-modules with the same $q$-dimension.
Hence we also obtain a well-defined  fusion tensor category associated to
$U_q\so_N\rtimes \Z_2$, with the usual restriction rules to $U_q\so_N$.

Recall the construction of spinors in the classical setting:  Consider a simple module of the Clifford algebra on $V=\C^N$.
It is well-known that for $N$ even we get an irreducible representation $S$ of
$Pin(N)$ which decomposes as a $Spin(N)$-module into a direct sum $S\cong S_+\oplus S_-$ of two non-isomorphic irreducible
representations. If $N$ is odd, we have two non-isomorphic
simple modules of the Clifford algebra, say $S_0$ and $S_1$, both of which restrict to
the same irreducible $Spin(N)$-module. We will just denote them (both) as $S$, consistent with the notation above.
For $N$ odd, we will also need the (reducible) $Pin(N)$-module $\tilde S=S_0\oplus S_1$ at some point.

The relationships between the spinor representations of $\U$ and $U_q\so_N$ are analogous to those of $Pin(N)$ and $Spin(N)$.
That is, for $N$ even, we have a $\U$-module $S$ which is irreducible and decomposes as $S\cong S_+\oplus S_-$
as a $U_q\so_N$-module.  For $N$ odd, there are two non-isomorphic $\U$-modules $S_0$ and $S_1$ which are isomorphic upon restriction to $U_q\so_N$ ($S_0$ and $S_1$ differ only on the $\Z_2$-action).

\subsection{Classical case} We first check some well-known identities in
the classical case, where $\U$ is replaced by $Pin(N)$ and $U_q'\so_n$ is replaced by $SO(n)$.
Most of these results have already more or less explicitly appeared,
as special cases of a more general approach, see \cite{Has}.

We consider the case where $Pin(N)$ representations are also $O(N)$
representations. We remark that
our symmetric bilinear form on the root lattice is normalized so that $\langle\beta,\beta\rangle=2$ for \emph{long} roots for all $N$,
for uniformity's sake. Recall (see e.g. \cite{Wy}) that simple $O(N)$
representations are labeled by the Young diagrams $\la$ for which
$\la_1'+\la_2'\leq N$ (here $\la_i'$ denotes the number of boxes in
the $i$-th column). The representations of the Lie algebra $\so_n$
for $n=2j$ are labeled by the dominant integral weights
$\mu=(\mu_i)_i$ such that $\mu_1\geq \mu_2\geq\ ...\ \mu_{j-1}\geq
|\mu_j|$, where either all $\mu_i$ are integers or all $\mu_i\equiv
1/2$ mod $\Z$. Then it is easy to check that the map
\begin{equation}\label{corresp}
\la\mapsto \bar\la,\quad{\rm where}\ \bar\la_i=N/2-\la'_{j+1-i}
\end{equation}
defines a bijection between the set of simple representations $V_\la$ of $O(N)$ for which $\la_1\leq n/2=j$ and the set of
simple $\so_n$ representations
$V_{\bar\la}$ for which $\bar\la_1\leq N/2$ and $N/2-\bar\la_i$ is an integer
for $1\leq i\leq n/2$.
Now consider the obvious action of $O(N)\times SO(n)$ on
$\C^N\otimes \C^n$. This induces commuting actions of $O(N)$ and $SO(n)$
via automorphisms on $Cliff(\C^N\otimes \C^n)$, and hence to projective
actions of these groups on a simple module $S_{Nn}$ of $Cliff(\C^N\otimes \C^n)$
i.e. proper actions of the corresponding covering groups, the spinor groups.

\begin{lemma}\label{comparison} (a) Let $n$ be even and let $S_{Nn}$ be
a simple module of the simple algebra $Cliff(\C^N\otimes \C^n)$.
Then $S_{Nn}$ decomposes as an $O(N)\times Spin(n)$ module as
$$S_{Nn}\cong \bigoplus_\la V_\la\otimes V_{\bar\la},$$
where $V_\la$ and $V_{\bar\la}$ are simple $O(N)$ and $Spin(n)$-modules and $\la$ runs through the set of Young diagrams
as in Eq \ref{corresp}.

(b) If both $N$ and $n$ are even, $S^{\otimes n}$ is isomorphic as a $Pin(N)\times Spin(n)$ module to the module $S_{Nn}$ in (a).
If $N$ is odd and $n$ is even,
$\tilde S^{\otimes n}$ is isomorphic as a $Pin(N)\times Spin(n)$ module
to the direct sum of $2^{n/2}$ copies of $S_{Nn}$ as in (a).

(c) Regardless of parity of $N$ and $n$, the irreducible representations of $Spin(n)$
in cases (a) and (b) are labeled by the dominant integrals weights $\mu$
satisfying $\mu_1\leq N/2$
and such that $\mu_i-N/2$ is an integer for $1\leq i\leq k$.
\end{lemma}

$Proof.$ It suffices to calculate the $Pin(N)\times Spin(n)$ characters of the various modules.
Let $n=2k$ and $\i=(i_1,\ldots,\ i_k)\in \Z_{\geq 0}^k$. We denote by
$\om(\i)$ the $Spin(n)$ weight given by the vector $(i_j-N/2)_j$.
Then we claim that the $Spin(n)$ character of a simple $Cliff(\C^N\otimes \C^n)$
module is given by
\begin{equation}\label{scharacter}
\chi^{S_{Nn}}=\sum_{i_1=0,\ldots,\ i_k=0}^N (\prod_{j=1}^k \chi_{i_j}) e^{\om(\i)},
\end{equation}
where $\chi_i$ is the $O(n)$ character for the $i$-th antisymmetrization
$\wedge^i V$ of the vector representation of $O(N)$, for $0\leq i\leq N$.
This can be seen as follows: As $Nn$ is even by assumption, we can describe
the character of the full spinor representation of $O(Nn)$ (which is
a simple $Cliff(Nn)$-module) by
$$(z_1z_2\cdots z_{Nn/2})^{-1/2}\sum_{j=0}^{Nn/2} e_j(z),$$
where $e_j(z)$ is the $j$-th elementary symmetric function in the variables
$z_1,\ldots ,z_{Nn/2}$. To view this as a character of $Spin(n)$
we replace the $z$-variables by variables $x_iy_j$, $1\leq i\leq n/2$,
$1\leq j\leq N$. We regard the result as a polynomial in the $x_i$ variables
over the ring of polynomials in the $y_i$ variables.
As every $x_i$ variable comes with all possible $y_j$
variables, and our formula is obviously symmetric in the $z$-variables,
and hence also in the $x$ and $y$ variables, a monomial in the
$x$-variables containing
the variable $x_i$ with the power $m_i$ must also have the factor
$e_{m_i}(y)$, the elementary symmetric function in the variables $y_1,\ldots, y_N$. Now it is well-known that the
elementary symmetric functions are the characters of the antisymmetrizations of the vector representation
which remain irreducible as $O(N)$-modules. This proves Eq \ref{scharacter}.

We can now prove statement (a) by induction with respect to inverse alphabetical order of the weights $\om(\i)$.
It is clear that the highest possible weight
occurring in Eq \ref{scharacter} is $\om=N\ve$. Then the coefficient of
$e^\om$ is equal to the trivial character, which proves (a) for $\la=0$.
The general claim follows by induction, using the formula
$$\prod\chi_{i_j}=\chi_\la+ lower\ characters,$$
where $i_j$ is a nonincreasing sequence of integers, $\la$ is the
Young diagram whose $j$-th column has exactly $i_j$ boxes, and lower
characters refers to a sum of simple $O(N)$ characters labeled by
Young diagrams smaller than $\la$ in alphabetical order.

To prove the corresponding formulas for the tensor product
representations, we check it first for $n=2$. Here for $N$ even, the
second tensor product of the spinor representation $S$ is a direct
sum of all possible antisymmetrizations of the vector representation
$\C^N$. For $N$ odd, we similarly get that $\tilde S^{\otimes 2}$
decomposes into the direct sum of two copies of the exterior algebra
of $\C^N$. It was shown in \cite{WSp} that the $i$-th
antisymmetrization in $S^{\otimes 2}$ (respectively in $\tilde S^{\otimes
2}$, where it appears with multiplicity 2) is an eigenspace of the
$\so_2$ generator $B_1$ with eigenvalue $N/2-i$. This proves that
the $\so_2$ character of $S^{\otimes 2}$ (respectively of $\tilde S^{\otimes
2}$) is given by Eq \ref{scharacter} for $N$ even (respectively by twice the
value of Eq \ref{scharacter} for $N$ odd). For $n=2k>2$, we write
$S^{\otimes n}=(S^{\otimes 2})^{\otimes k}$ and observe that the
$i$-th factor $S^{\otimes 2}$ gives us the eigenspaces of
$B_{2i-1}$, to which we can apply the same arguments as before.
Comparing with Eq \ref{scharacter} (with the $\chi_{i_j}$ evaluated
at the identity element) we see that the $SO(n)$ character of
$S^{\otimes n}$ is the same as the one for $S_{Nn}$ for $N$ even,
and the $SO(n)$ character of $\tilde S^{\otimes n}$ is $2^{n/2}$
times the character of $S_{Nn}$ for $N$ odd. From this follow
statements (b) and (c) (for $n$ even). For $n$ odd, the
corresponding statements follow from the results for $n+1$ from the
restriction rules of representations of $\so_{n+1}$.\qed

\subsection{Quantum and fusion cases} By the main result of \cite{WSp}, we have commuting actions of
$\U=U_q\so_N\rtimes \Z/2$ and $U'_q\so_n$ on $S^{\otimes n}$ (for $N$ even) and $\tilde S^{\otimes n}$ for $N$ odd.
Not surprisingly, the decomposition in the Lemma \ref{comparison} carries over to this setting if $q$ is not a root of unity.
If $q$ is a primitive $2\ell-$th root of unity, we have a similar relationship in the corresponding ribbon fusion category $O(N)_r$ where $r=\ell+2-N$ is the level.
This is the quotient category of the (ribbon) category of tilting modules in  $\U=U_q\so_N\rtimes \Z_2$ by negligible morphisms.
Adopting the notation from the affine Lie algebra literature, we denote this category by $O(N)_r$ where $r=\ell+2-N$.
In the case $N$ is odd, we have $O(N)_r\cong SO(N)_r\boxtimes \Rep(\Z_2)$, whereas in the case $N$ is
even $O(N)_r$ is the $\Z_2$-equivariantization of $SO(N)_r$.  The simple objects in $O(N)_r$
corresponding to $O(N)$-representations are labeled by Young diagrams $\la$
satisfying $\la_1'+\la_2'\leq N$ and $\la_1+\la_2\leq \ell +2-N$ and the additional Young diagram $\la=[\ell-N+2,1^{N}]$.
The objects with half-integer spin can be described by similar inequalities.
A more explicit description is given below in the case $r=2$.
We will again denote the images of the corresponding tilting modules in $\U$ by $S$ (respectively $\tilde S$) in the fusion category $O(N)_r$.
We have the following results, most of which were already proved in \cite{WSp}:

\begin{theorem}\label{duality} (a) Let $n$ be even. Then we can define an action of $\U\times U'_q\so_n$ on  $S^{\otimes n}$ for $N$ even (respectively
$\tilde S^{\otimes n}$ for $N$ odd) whose decomposition into irreducibles is the same as in the classical case,
if $q$ is not a root of unity.

(b) If $q$ is a primitive $2\ell$-th root of unity, then the objects
$S^{\otimes n}$ for $N$ even (respectively $\tilde S^{\otimes n}$ for $N$ odd) decompose in $O(N)_{\ell-N+2}$ as a direct sum
$\bigoplus_\la V_\la\otimes V_{\bar\la}$. Here now
$V_\la$ ranges over the objects as in the classical case, subject to
the additional condition $\la_1+\la_2\leq \ell +2-N$, and the additional diagram $\la=[\ell-N+2,1^{N}]$, and
$V_{\bar\la}$ is the (via \ref{corresp}) corresponding  $U'_q\so_n$
module with highest weight $\bar\la$.
\end{theorem}

$Proof.$  Part (a) follows from Lemma \ref{comparison}, using the
explicit representations in \cite{WSp} and the fact that for $q$ not
a root of unity the representation theory of Drinfeld-Jimbo quantum
groups is essentially the same as for the corresponding Lie algebra.
For part (b) we just use the fact that tensor powers of $S$ and
$\tilde S$ can be written as a direct sum of indecomposable tilting
modules; the objects in the fusion category are obtained by taking
the quotient module by the tensor ideal generated by those tilting
modules which have $q$-dimension equal to 0. The representations of
$U'_q\so_n$ into these tensor powers are still well-defined at a
root of unity, and they factor over the fusion quotient. As these
$U'_q\so_n$ modules usually have smaller dimensions at a root of
unity than in the generic case, we still need to check that they
have the same highest weight vector. But this follows from the
restriction rule: restricting the action to $\so_{n-1}$, the highest
weight vector is again a highest weight vector in an
$\so_{n-1}$-module which also exists in the fusion category. The
explicit combinatorics can be checked either directly by using
Gelfand-Tseitlin bases for the orthogonal case (see e.g. \cite{GK}),
or by using the tensor product rules for spinor representations (see
e.g. \cite{WSp}) via the correspondence \ref{corresp}.\qed
\medskip

We use the notation $\ve=(1/2,1/2,\ldots,1/2)\in
\R^j$ and $\e_i$ for the $i$-th standard basis vector of $\R^j$. We associate
these vectors with weights of $\so_n$ for $n=2j$ or $n=2j+1$ in the usual way.
Let $\rho$
be half the sum of the positive roots of $\so_n$, and let
$q^{2\rho}$ be the operator on a finite dimensional $\U$-module
defined by $q^{2\rho}v_\mu=q^{(2\rho, \mu)}v_\mu$ for a weight
vector $v_\mu$ of weight $\mu$. We define, as usual, the
$q$-dimension of a $\U$-module $V$ by $\dim_q V=Tr(q^{2\rho})$. As
we have commuting actions of $\U$ and $U'_q\so_n$ on $S^{\otimes n}$
(respectively $\tilde S^{\otimes n}$), we can define the virtual $U'_q\so_n$
character $\chi^\rho_n$ by
$$\chi^\rho_n(u)=Tr(uq^{2\rho}),$$
where $u$ is in the Cartan algebra of $U'_q\so_n$, and $Tr$ is the
usual trace of $S^{\otimes n}$ (respectively $\tilde S^{\otimes n}$). The
following lemma follows from the multiplicativity of the trace for
tensor factors, using a similar argument as in the proof of Lemma
\ref{comparison}.

\begin{lemma}\label{rhocharacter} If $N$ is even
then the character $\chi^\rho_n$ is uniquely determined by
$$\chi^\rho_n(\prod
B_{2i-1}^{e_{2i-1}})=\prod \chi^\rho_2(B_{2i-1}^{e_{2i-1}}),$$ and
$\chi_2^\rho(B_1^e)=\sum^N_{j=1} \dim_qV_{[1^j]}[N/2-j]^e$. If
$N$ is odd, the same formulas hold, except that we have to add a
factor 2 on the right hand side of the formula for each
$\chi_2^\rho(B_{2i-1}^{e_i})$, $1\leq i<N/2$.
\end{lemma}

\subsection{Weakly Integral Cases}

In the rest of this paper we will mostly focus on the case $q=e^{\pi \i/N}$ corresponding to $O(N)_2$.

The special cases $O(N)_2$ correspond to the quotient by negligible morphisms of the categories of tilting $\U$- modules
for $q$ a $2N$th root of unity.  These $O(N)_2$ are weakly integral unitary ribbon fusion categories,
i.e. $(\dim_q V)^2\in\Z$ for simple objects $V$.

The related categories $SO(N)_2$ (see e.g.
\cite{NR}) obtained from $U_q\so_N$ at $q=e^{\pi \i/N}$ are also weakly integral modular categories and have simple objects
labeled by highest weights for $\so_N$. We will describe these categories in some detail.

 Setting
$N=2k+1$ for $N$ odd and $N=2k$ for $N$ even, we denote the fundamental weights for $\so_N$ by
$\Lambda_1,\ldots,\Lambda_k$. No confusion should arise as we
deal with $N$ even and $N$ odd separately.  For later use we define for $0\leq
j\leq k$ the highest weight $\gamma_j=(1,\ldots,1,0,\ldots,0)$ with the first $j$ entries equal to $1$.

For $N$ odd $\Lambda_k=(1/2,\ldots,1/2)$ labels
the simple object $S$ associated with the fundamental spin representation for
$\so_N$ and $\Lambda_j=(1,\ldots,1,0,\ldots,0)$ for $1\leq j\leq k-1$.

For $N$ even the two fundamental spin objects $S_{\pm}$ are
labeled by $\Lambda_k=(1/2,\ldots,1/2)$ and $\Lambda_{k-1}=(1/2,\ldots,1/2,-1/2)$, while
$\Lambda_j=(1,\ldots,1,0\ldots,0)$ for $1\leq j\leq k-2$.

\subsubsection{$N$ odd}\label{ss:nodd}
The fusion category $SO(N)_2$ for $N$ odd has two simple (self-dual) objects
$S=V_{\Lambda_k}$ and $S^\prime=V_{\Lambda_k+\Lambda_1}$ of
dimension $\sqrt{N}$, 2 simple objects $\1=V_{\gamma_0}$ and $V_{2\Lambda_1}$ of
dimension $1$, and
$\frac{N}{2}$ simple objects $V_{\gamma_s}$ of dimension $2$ where $1\leq s\leq
\frac{N-1}{2}$.  Thus, for $N$ odd, the rank of $SO(N)_2$ is
$\frac{N-1}{2}+4$ and the categorical dimension is $4N$.

As we have noted above, for $N$ odd $O(N)_2\cong SO(N)_2\boxtimes \Rep(\Z_2)$ as ribbon fusion categories,
so that the structure of $O(N)_2$ is easily determined from that of $SO(N)_2$.  Here $\Rep(\Z_2)$ is regarded as the ribbon category with trivial twists and symmetric braiding.  We will denote the two objects in $\Rep(\Z_2)$ be $\one$ and $-\one$ where $-\one$ corresonds to the non-trivial representation of $\Z_2$.  In particular we have a $\Z_2$ grading of $O(N)_2$ with components corresponding to $(V,\pm\one)$.  In this notation we have $S_0=(S,\one)$ and $S_1=(S,-\one)$.  For example we have
$\tilde{S}^{\ot 2}=[(S,\one)\oplus (S,-\one)]^{\ot 2}\cong 2[(S^{\ot 2},\one)\oplus (S^{\ot 2},-\one)]$.  Moreover, the (forgetful) functor $F:O(N)_2\rightarrow SO(N)_2$ by $F(V,\pm\one)\rightarrow V$ is obviously faithful and is braided since the braiding on $\Rep(\Z_2)$ is symmetric.

\subsubsection{$N$ even}
For $N$ even the fusion category $SO(N)_2$ has $4$ simple objects $S_{\pm}$
(labeled by $\Lambda_k$ and $\Lambda_{k-1}$) and
$S_{\pm}^\prime$ (labeled by $\Lambda_k+\Lambda_1$ and
$\Lambda_{k-1}+\Lambda_1$) of dimension $\sqrt{N/2}$, $4$ simple objects
$\1,V_{2\Lambda_1},V_{2\Lambda_k}$, and $V_{2\Lambda_{k-1}}$ of dimension $1$
and
$\frac{N}{2}-1$ simple objects $V_{\gamma_s}$ of dimension $2$ where $1\leq
s\leq
\frac{N}{2}-1$.  Thus, for $N$ even, the rank of $SO(N)_2$ is $\frac{N}{2}+7$ and the
categorical dimension is $4N$.

 The simple
objects in $O(N)_2$ are the images (under purification) of the simple $U_q\so_N\rtimes
\Z_2$-tilting modules with non-zero $q$-dimension.  Using \cite[Section 3.4]{WSp} we find that the simple objects
in $O(N)_2$ are: $S$ and $S^\prime$ of dimension $2\sqrt{N/2}$;
$\1,V_{[2]},V_{[1^N]}$ and $V_{[1^{N-1},1]}$ of dimension $1$; and $V_{[1^s]}$
of dimension $2$ with $1\leq s\leq N-1$.  The restriction map
$\Rep(U_q\so_N\rtimes \Z_2)\rightarrow \Rep(U_q\so_N)$ induces a braided
tensor functor $F:O(N)_2\rightarrow SO(N)_2$ with images:
\begin{eqnarray*}
 F(S)&=&S_+\oplus S_-\\
F(S^\prime)&=& S_+^\prime\oplus S_-^\prime\\
F(V_{[2]})&=&F(V_{[1^{N-1},1]})=V_{2\Lambda_1}\\
F(V_{[1^N]})&=&F(\1)=\1\\
F(V_{[1^s]})&=&F(V_{[1^{N-s}]})=V_{\gamma_s}, \quad 1\leq s\leq k-1
\end{eqnarray*}

Observe that the objects $S$ and $S^\prime$ in $O(N)_2$ are self-dual, although $S_{\pm}$ are not.

We now can proof the following corollary to Theorem \ref{duality}(b).

\begin{corollary}\label{fusion2}
Let $q$ be a primitive $2N$-th root of unity. Then the representations $\Phi$ of $U'_q\so_n$ into the $O(N)_2$
centralizer algebras $\End(S^{\otimes n})$ for $N$ even (respectively $\End(\tilde S^{\otimes n})$ for $N$ odd)
are labelled by the weights $N\ve$ and $N\ve-\e_j$ if $n=2j+1$ is odd,
and by the weights $N\ve -r\e_j$,
$0\leq r\leq N$, $N\ve-\e_{j-1}-\e_j$ and $N\ve-\e_{j-1}-(N-1)\e_j$
if $n=2j$ is even.  For $N$ even, $\Phi$ is surjective.
\end{corollary}

$Proof.$ This follows from Theorem \ref{duality}(b) and from the restriction rules for
representations of $U'_q\so_n$ (see \cite[Lemma 4.2 and Proposition 4.3]{WSp}) and tensor product rules of $O(N)_2$.  The surjectivity for $N$ even follows from a dimension count (simply compute the Bratteli diagram for the object $S$).\qed

\begin{remark}\label{nosurjrem} \begin{enumerate}
\item[(a)] One can check that already for $n=2$
the map $\Phi$ is not surjective for $N$ odd.  It follows from the
explicit representation of $U'_q\so_n$ in \cite{WSp} that the image
of $B_1$ in $\tilde{S}^{\ot 2}=(S_0\oplus S_1)^{\otimes 2}=$ permutes $S_1^{\otimes 2}$
and $S_0^{\otimes 2}$, and similarly for the mixed terms.  In particular,
 $\Phi(B_1)$ does not commute with the projections $p\in\Hom(\tilde{S}^{\ot 2},V)$ where $V$ is a simple subobject of $\tilde{S}^{\ot 2}$.  However, $\Phi(B_1^2)$ leaves $S_i^{\ot 2}$, $S_0\ot S_1$ and $S_1\ot S_0$ invariant (see \cite{WSp}), and more generally $\Phi(B_i^2)$ leaves invariant any tensor product of any number of copies of $S_0$ and $S_i$.
 \item[(b)] Using the notation $(V,\pm\one)$ for objects in the two components of the $\Z_2$-grading on $O(N)_2$, we have $\tilde{S}^{\ot n}\cong 2^{n-1}(S^{\ot n},\one)\oplus 2^{n-1}(S^{\ot n},-\one)$.  The (faithful braided tensor) functor $F:O(N)_2\rightarrow SO(N)_2$ induces an algebra homomorphism $\Xi:\End(\tilde{S}^{\ot n})\rightarrow\End((2S)^{\ot n})=\End(2^{n-1}S^{\ot n}\oplus 2^{n-1}S^{\ot n})$.  The image of $\Xi$ lies in the diagonal:
 $\End(2^{n-1}S^{\ot n})\times\End(2^{n-1}S^{\ot n})$, since $\Hom((S^{\ot n},\one),(S^{\ot n},-\one))=0$.  Moreover, as $\Phi(B_i^2)\in\End(\tilde{S}^{\ot n})$ leaves invariant each of the $2^{n-1}$ copies of $(S^{\ot n},\pm\one)$, we see that $\Xi(\Phi(B_i^2))$ lies in the diagonal $\prod_{i=1}^{2^n}\End(S^{\ot n})$.  Since the latter algebra is isomorphic to $\End(S^{\ot n})$ we see that $\Xi\circ\Phi(\langle 1,B_1^2,\ldots,B_{n-1}\rangle)$ is isomorphic to a subalgebra of $\End(S^{\ot n})$.
                                \end{enumerate}

\end{remark}

\subsection{$\B_n$ representations on $\End(S^{\ot n})$}\label{quantugroupeigs}

Denote by $\gamma_S:\B_n\rightarrow \Aut(S^{\ot n})$ the representations of the braid group associated with the object $S$ in
$SO(N)_2$ for $N$ odd or $O(N)_2$ for $N$ even.  Explicitly, $\gamma_S$ is defined
on generators by $\sigma_i\rightarrow Id_S^{\ot(i-1)}\ot c_{S,S}\ot Id_S^{\ot (n-i-1)}$.

For later use we compute the eigenvalues for the braiding operator $c_{S,S}$
for $SO(N)_2$ when $N$ is odd and $O(N)_2$ for $N$ even.

\begin{remark}
In the subsection only, we let $x=e^{\pi \i/(2N)}$, and let $\langle\;,\;\rangle$ be the symmetric bilinear
form on the weight lattice normalized so that $\langle\al,\al\rangle=2$ for
short roots.  This is to conform with the standard results, and is only different for $N$ odd.
\end{remark}

For $N=2k+1$ odd,  we have
$$S^{\ot
2}\cong\bigoplus_{j=0}^{k} V_{\gamma_j}.$$
The eigenvalues of $c_{S,S}$ are easily computed, and we record them in:
\begin{lemma}
Let $N=2k+1$ be odd. Up to an overall factor depending only on $N$, the
eigenvalue of $c_{S,S}$ on the projection onto the simple object $V_{\gamma_s}$
is
\begin{equation}\label{Noddeigs}
 \Psi(N,s):=i^{(k-s)^2}e^{-\pi \i s^2/(2N)}.
\end{equation}
\end{lemma}
$Proof.$
It follows from
Reshetikhin's formulas (see e.g.
\cite[Corollary 2.22]{LR})
that, up to an overall factor, $c_{S,S}$ acts on the projection onto
$V_{\lambda}$ by the scalar
$\varsigma(\lambda)x^{\frac{c_{\lambda}}{2}}$ where
$c_\lambda=\langle \lambda+2\rho,\lambda\rangle$ for any weight $\lambda$ and
the sign $\varsigma(\lambda)=1$ if the corresponding $\so_N$ representation
appears in the symmetric tensor square of the fundamental spin representation
and $-1$ otherwise. Observe that here $\langle\;,\;\rangle$ is twice the usual
Euclidean inner product and $2\rho=(2k-1,\ldots,1)$.
We compute $c_{\gamma_s}=2(Ns-s^2)$ and note that
$$\varsigma(\gamma_s)=\begin{cases}
                      -1 & (k-s)\equiv 1,2\pmod{4}\\
1 & (k-s)\equiv 0,3\pmod{4}
                     \end{cases},$$ from which the result follows.\qed

In the case $N=2k$ is even we have:
$$S^{\ot 2}=\bigoplus_{s=0}^N V_{[1^s]}$$
and the eigenvalues of $c_{S,S}$ are given in:
\begin{lemma}\label{Neveneigs}
Let $N=2k$ be even. Up to an
overall factor depending only on $N$, the eigenvalue of the
$O(N)_2$ braiding operator $c_{S,S}$  on the projection onto the simple object
labeled by $[1^s]$ is:
$\eta(s)f(s)$, where $\eta(s)=e^{(N-2s)(N-2s+2)\pi \i/8}$ and
$f(s)=i^se^{-\pi \i s^2/(2N)}$.
\end{lemma}
$Proof.$
Since the functor $F:O(N)_2\rightarrow SO(N)_2$ is a braided tensor functor we
can compute the eigenvalues of $c_{S,S}$ from $F(c_{S,S})$.  Up to signs these
are just the eigenvalues of $c_{S_\pm,S_\pm}$ and the
square roots of the eigenvalues of $c_{S_+,S_-}c_{S_-,S_+}$.  These can be
computed up to an overall factor  using Drinfeld's quantum Casimir \cite{Dr} (since $\langle
\Lambda_k+2\rho,\Lambda_k\rangle=\langle
\Lambda_{k-1}+2\rho,\Lambda_{k-1}\rangle$) as $q^\frac{c_\lambda}{2}$ with
$q=e^{\pi \i /N}$ for any $V_\lambda\in F(S^{\ot 2})$. Up to signs, the
eigenvalues corresponding to $V_{[1^s]}$ and $V_{[1^{N-s}]}$ are (both)
$q^{\frac{c_{\gamma_s}}{2}}$ for $0\leq s\leq N/2$.  We compute
$c_{\gamma_s}=\langle \gamma_s+2\rho,\gamma_s\rangle=Ns-s^2$ and set
$f(s)=q^{\frac{Ns-s^2}{2}}=i^se^{-\pi \i  s^2/(2N)}$.  Observe that $f(N-s)=f(s)$
so that $c_{S,S}$ has eigenvalue $\eta(s)f(s)$ on the projection onto
$V_{[1^s]}$ for all $0\leq s\leq N$, where $\eta(s)$ is a sign.

By continuity, it is enough to determine the signs for the classical case $q=1$ for which the braiding is symmetric.
One way to do this goes by induction on the dimension $N$, for $N$ even (a similar argument also works for
the slightly easier case $N$ odd). One first observes that
for $N=4$ the signs are given by $\eta(0)=\eta(1)=\eta(4)=-1$ and $\eta(2)=\eta(3)$,
using the fact that $Spin(4)\cong SU(2)\times SU(2)$.

The crucial observation now is that the sign
for the representations $V_{[1^{N/2-s}]}\subset S_N^{\otimes 2}$
 are the same as the ones for the representations
$V_{[1^{N/2-s-1}]}\subset S_{N-2}^{\otimes 2}$, for $0\leq |s|<N/2$;
here $S_{2k}$ is the spinor representation in connection with $O(2k)$.
This follows from the fact that $S_N$ decomposes as a $Pin(N-2)$ module
into the direct sum of two modules isomorphic to $S_{N-2}$, see e.g. the discussion in \cite{WSp}, Lemma 2.1.
Using the eigenspace decomposition of the permutation $R_S\in\End(S^{\otimes 2})$,
we obtain for the normalized trace $tr$ on $\End(S^{\otimes 2})$
\begin{equation}\label{traceformula}
\frac{1}{2^{N/2}}=tr(R_S)=\frac{1}{2^{N}}\sum_{s=0}^N \eta(N/2-s)\dim\ V_{[1^{N/2-s}]}.
\end{equation}
We remark that a similar formula also holds for the odd-dimensional case
$Spin(N+1)$, where now the summation only goes until $s=N/2$ and
we have the antisymmetrizations of the $(N+1)$-dimensional vector representations on the right hand side.
By induction assumption, $\eta(N/2-s)$ is known for $s<N/2$, and
$\dim V_{[1^{N/2-s}]}$ is equal to $\binom{N}{N/2-s}$.
In the odd-dimensional case, we can now easily calculate the missing sign
 $\eta(0)$ from Eq \ref{traceformula}, as adjusted for the odd-dimensional case. To calculate
the two remaining signs in the even-dimensional case, we consider $Pin(N)$ as a subgroup of
$Spin(N+1)$, which acts irreducibly
via its spinor representation on the same
vector space $S$; in particular, we can also identify the trivial
subrepresentation in $S^{\otimes 2}$ for both groups, which hence has
the same sign $\eta(0)$ for the permutation $R_S$ at $q=1$. One now calculates
$\eta(N)$ from  Eq. \ref{traceformula}.
It is now easy to check that the signs can be given by the formula $\eta(s)=e^{\frac{(N-2s)(N-2s+2)\pi \i }{8}}$.

\qed

\section{Representation theory of $U_q'\so_n$}\label{section: rep. theory}

We review and (re)prove certain results of the representation theory
of $U_q'\so_n$. Many of these results have already appeared in one
form or another in work of Klimyk and his coauthors, see e.g.
\cite{GK}, \cite{IK}. However, in our case, we need these results
for roots of unity where the situation is more complicated.
Hence we have decided to give our own, quite different proofs by
mimicking a Verma module construction. We will do this here only for
what is called the classical series in \cite{IK}, i.e. for
representations which are deformations of representations of
$U\so_n$, and those only for $n\leq 5$. It is planned to give a more
complete study of these representations in a separate paper \cite{Wrep}.

\subsection{Definitions} We identify roots and weights of $U_q'\so_n$ with
vectors in $\R^k$, where $k=n/2$ or $(n-1)/2$ depending on the parity of $n$,
as usual. So if $\ep_i$ is the $i$-th standard unit vector for $\R^k$, the roots are given
by $\pm \ep_i\pm \ep_j$, $1\leq i<j\leq k$, and, if $n=2k+1$ is odd,
also by $\pm \ep_i$, $1\leq i\leq k$. Here the analog of the Cartan subalgebra is the algebra $\h$ generated by
$B_1, B_3,\ldots, B_{2k-1}$ for $n=2k$ or $n=2k+1$. A vector $v$ in a $U_q'\so_n$-module is said to have
weight $\la$ if $B_{2i-1}v=[\la_i]v$
for all $B_{2i-1}\in \h$; we shall often identify $\la$ with the vector
$(\la_i)$. As usual, $[n]=(q^n-q^{-n})/(q-q^{-1})$.
Let us first recall the following theorem, which has been proved in \cite{IK};
it also follows from the results in \cite{WSp}, as quoted in Theorem \ref{duality}.

\begin{theorem}\label{genericrep}
Let $\la$ be a dominant integral weight, and let $q$ be generic. Then there exists a finite dimensional simple
$U_q'\so_n$ module $V_\la$ with highest weight $\la$ and the same weight multiplicities as for the corresponding $U\so_n$ module.
\end{theorem}

\begin{lemma}\label{weightshift}
Let $v$ be a vector in a $U_q'\so_n$ module with weight $\mu$. Then

(a)\ $(B_{2i-1}-[\mu_i+1])(B_{2i-1}-[\mu_i-1])B_{2i}v=0$,

(b)\ $(B_{2i+1}-[\mu_{i+1}\pm 1])(B_{2i-1}-[\mu_i+1])B_{2i}v$ has weight $\mu-(\ep_i\pm \ep_{i+1})$, if it is nonzero.

\noindent In particular, we can write $B_{2i}v$ as a sum of two eigenvectors of $B_{2i-1}$ (if $[\mu_i+1]\neq [\mu_i-1]$),
and we can write $(B_{2i-1}-[\mu_i+1])B_{2i}v$ as a sum of two weight vectors (if $[\mu_{i+1}+1]\neq [\mu_{i+1}-1]$).
\end{lemma}

$Proof.$ These are straightforward calculations. E.g. for (a) we have
\begin{align*}
B_{2i-1}^2B_{2i}v\ &=\ ([2]B_{2i-1}B_{2i}B_{2i-1}-B_{2i}B_{2i-1}^2+B_{2i})v\cr
&=\ [2][\mu_i]B_{2i-1}B_{2i}v-([\mu_i]^2-1)B_{2i}v.
\end{align*}
We now get the claimed factorization in (a) using the identities
$[2][\mu_i]=[\mu_i+1]+[\mu_i-1]$ and $[\mu_i]^2-1=[\mu_i+1][\mu_i-1]$.
For part (b)  observe that a similar calculation also holds
with $i$ replaced by $i+1$. The claim follows from this.\qed
\medskip

 For a given weight $\la$ we define the left ideal
\begin{equation}\label{Ilambdadef}
I_\la=U'_q\so_n\langle (B_{2i-1}-[\la_i]1), (B_{2i-1}B_{2i}-[\la_i-1]B_{2i})\rangle
\end{equation}
for all values of $i$ for which the indices $2i-1$ and $2i$ are between
(including) 1 and $(n-1)$. Observe that one can show as in Lemma \ref{weightshift} that now $B_{2i}$ is an eigenvector of
$B_{2i-1}$ with eigenvalue $[\la_i-1]$ mod $I_\la$. Moreover, if $[\la_{i+1}+1]\neq [\la_{i+1}-1]$, we can write $B_{2i}$
as a linear combination of
the two vectors $(B_{2i+1}-[\la_{i+1}\pm 1])B_{2i}$ of
weights $\la-(\ep_i\pm\ep_{i+1})$  mod $I_\la$;
observe that these are weights of the form $\la-\al$ with $\al$ a positive root of $U'_qso_n$.

\subsection{Spanning property} It has already been observed in
\cite{GK} that a PBW type theorem holds for the algebra $U'_q\so_n$, using
its embedding into the quantum group $U_q\sl_n$. One can also prove the existence of an analogue of a Verma module.
This, and more results, are planned to appear in a separate paper \cite{Wrep} by the second named author.
For this paper,
we will only give (or outline) \emph{ad hoc} proofs for the special cases needed for our purpose.

\def\Vlt{\tilde{M}}
\def\S{{\mathcal S}}
\begin{lemma}\label{lemrelations}
Let $\la$ be a weight of $U'_q\so_n$ for $n\leq 5$. Then the Verma module $\Vl=U'_q\so_n/I_\la$ is spanned by the ordered
products of the form
$$B_2^{e_1}(B_3B_2)^{e_2}(B_4B_3B_2)^{e_3}B_4^{e_4}v_0,$$
where the $e_i$ are nonnegative integers which, for $n<5$, are equal to 0 for those factors which are not in $U'_q\so_n$,
and where $v_0 \equiv 1$ mod $I_\la$ is the highest weight vector.
\end{lemma}

$Proof.$ The proof can be done via elementary, albeit somewhat tedious, calculations.
A more general result will be proved in \cite{Wrep}. We give a fairly detailed
outline for a proof of this lemma for the skeptical reader as follows:

For $n=5$, the idea is to move
the generators $B_4$ as far to the right as possible. To make this mathematically precise,
we define an order on words in the generators $B_i$
first by the length of the word, and then by reversed alphabetical order
e.g. $B_4^2<B_3B_4<B_4B_3$ etc.
We first prove that the claim holds if we only apply generators
$B_i$, $2\leq i\leq 4$, to the highest weight vector.
As a first step
one shows that any vector generated this way is a linear combination
of vectors of the form $w(B_4B_3B_2)^{e_3}B_4^{e_4}v_0$, with the word $w\in \langle B_2, B_3\rangle$.
This follows by moving generators $B_4$ as far to the right as possible,
using the relation
$$B_4(B_4B_3B_2)\ =\ [2](B_4B_3B_2)B_4-B_3B_2B_4^2+B_3B_2.$$
It is not hard to show that one can express $B_3B_4^jv_0$ as a linear
combination of vectors $B_4^iv_\la$, see Lemma \ref{sothreeVerma} for details.
Moreover, we also have the relation
$$B_3(B_4B_3B_2)\ =\ (B_4B_3B_2)B_3+[B_3,B_2B_3B_4].$$
Using it, not only can we prove our claim, but we can also show that $w$ may be assumed
to end with a $B_2$, by induction on $e_4$ and $e_3$.
It is now an easy induction on the number of $B_3$s in $w$ to
prove that it can be expressed as a linear combination of words of
the form $B_2^{e_1}(B_3B_2)^{e_2}$ by moving the $B_3$s as far to the
right as possible (taking into account that a $B_3$ on the right end
of $w$ will be absorbed, as just mentioned).
To finish the proof for $n=5$, it
suffices to show that multiplying any of the words as in the
statement by $B_1$ again results in a linear combination of words
without a $B_1$; this follows by a similar induction on the order of
the words. The claims for $n=4$ and $n=3$ are proved similarly, with
the proofs being much easier. \qed

\begin{corollary}\label{weights} A weight appears in the highest weight module $N_\la$ for $U'_q\so_n$
with at most the multiplicity as in the Verma module $\Vl$ for the classical case $U\so_n$ at $q=1$, for $n\leq 5$.
\end{corollary}

$Proof.$ We give an outline of the proof for the most difficult case $n=5$.
As $N_\la$ is a quotient of $\Vl$, it suffices to prove the statement for the latter module.
It is standard to check that the elements $B_2,B_3B_2, B_4$ and $B_4B_3B_2$
form a basis of $(\so_5 +I_\la)/I_\la$ for $q=1$. Hence their ordered polynomials form a basis for $\Vl=U\so_n/I_\la$.

Let us consider the subspaces $\Vl(f_1,f_2)$ spanned by all the monomials in the generators
with at most $f_1$ and $f_2$ factors equal to $B_2$ and $B_4$ respectively. It follows from Lemma \ref{lemrelations}
and its proof
that any such element can be written as a linear combination of words which also contain $\leq f_1$ factors equal to $B_1$
and $\leq f_2$ factors equal to $B_{3}$. Hence this space is a module of the Cartan algebra generated by $B_1$ and $B_3$.
By Lemma \ref{weightshift},
the zeroes of the characteristic polynomial of $B_{2i-1}$ acting on
$\Vl(f_1,f_2)$ can only be of the form $[\la_i-j]$ for some integer $j$.
Specializing at $q=1$ gives us the estimates on the multiplicities of the zeroes (In fact, with a little more effort,
one could show that our basis for $q=1$ extends to a basis for general $q$, which proves equality for the multiplicities).
The general claim now follows by letting $f_1$ and $f_2$ go to infinity. \qed

\begin{remark} Having an analog of Verma modules, one can show that there exists a unique simple $U'_q\so_n$ highest weight module
with given highest weight, by the usual standard arguments, in the
generic case. Unfortunately, we will need this at roots of unity.
Results for the usual quantum groups at roots of unity would suggest
that there could be many nonisomorphic simple modules with the same
highest weight, see \cite{dC} and the papers quoted therein. This leads to
the consideration of certain invariant forms.
\end{remark}

\subsection{Invariant forms}
We call a sesquilinear form $(\ ,\ )$ on a $U'_q\so_n$ module $M$ \emph{invariant} if
$(B_iv,w)=(v,B_iw)$ for all $v,w\in M$ and $1\leq i<n$. A $U'_q\so_n$ module $M$ is called \emph{unitarizable}
if it admits a positive definite invariant form.

In the following, we will denote a highest weight module with highest weight $\la$ by $N_\la$.
If $q$ is a root of unity, the action of the operators $B_i$ on $N_\la$ may no longer be diagonalizable.
However, we only have finitely many (generalized) weight spaces.
For a weight $\mu$ we let $N_\la[\mu]$ be the generalized weight space of $N_\la$, i.e.
the set of all vectors $v$ such that $(B_{2i-1}-[\mu_i]1)^kv=0$ for sufficiently large $k$.
Finally, if $q$ is a primitive $2\ell$-th root of unity, with $\ell\geq n$,
we say that $\la$ is a {\it restricted dominant weight} for $U'_q\so_n$ if $\la_1\leq \ell/2$.

\begin{lemma}\label{bilinunique} Let $\la$ be a dominant integral weight with corresponding highest weight module $N_\la$ and highest weight vector $v_\la$.

(a) For $q$ not a root of unity, there is at most one invariant bilinear form $(\ ,\ )$ on $N_\la$, up to scalar multiples.

(b) Let now $q$ be arbitrary, and suppose $N_\la$ admits an invariant bilinear form $(\ ,\ )$. For $a=B_{i_1}B_{i_2}\cdots B_{i_k}$, set
$a^t=B_{i_k}\cdots B_{i_2}B_{i_1}$. Then the value of $(av_\la, bv_\la)$ is uniquely determined by $(v_\la,v_\la)$
whenever $a^tbv_\la$ can be written as a linear combination of
generalized weight vectors such that the $N_\la[\la]$ component is a multiple of $v_\la$.

\end{lemma}

$Proof.$ Part (a) follows from a standard argument, which we omit.  It follows from invariance that
$$(a_1v_\la,a_2v_\la)\ =\ (v_\la, a_1^ta_2v_\la).$$
If $q$ is not a root of unity, all the weight spaces are mutually orthogonal with respect to an invariant bilinear form. Hence the value of $(a_1v_\la,a_2v_\la)$ is given by the scalar of $v_\la$ in the expansion
of $a_1^ta_2v_\la$ as a linear combination of weight vectors,  times $(v_\la,v_\la)$. Part (b) is proved the same way. \qed

\begin{remark}\label{strategy} The strategy now will be to show that for certain dominant weights $\la$
there exists at most one
unitarizable simple module with highest weight $\la$. The idea is to
show that, loosely speaking, any additional vectors in the weight
space of $\la$ in $\Vl$ already have to be in the annihilator ideal
of a positive semidefinite form on $\Vl$.
\end{remark}

\subsection{$\so_3$} We now give a detailed classification of certain $U'_q\so_3$ modules
as these results will be used later for $q$ a root of unity.

\begin{lemma}\label{sothreeVerma} Suppose $q$ is not a root of unity, and define $v_0\in\Vl$ by $v_0=1$ mod $I_\la$ for the $U'\so_3$ weight $\la\in\R$. Then the set $\{B_2^jv_0,\ j\geq 0\}$ forms a basis of $\Vl$.
Moreover, $\Vl$ also has a basis of weight vectors $v_j$ with weight $[\la-j]$, $j=0,1,\ldots$  defined inductively by  $v_1=B_2v_0$ and
$$v_{i+1}=B_2v_i-\al_{i-1,i}v_{i-1},\quad {\rm for}\ i>1,$$
where
$$\al_{i-1,i}=\frac{[i][2\la -i+1]}{(q^{\la-i}+q^{i-\la})(q^{\la-i+1}+q^{i-\la-1})}.$$
In particular, if $\la$ is a half-integer, there exists a unique
simple module with highest weight $\la$ whose dimension is $2\la +
1$, and on which both $B_1$ and $B_2$ act with the same set of
eigenvalues $\{ [\la-j], 0\leq j\leq 2\la\}$. Finally, there is at
most one invariant form $(\ ,\ )$
on $\Vl$, up to scalar multiples. It is completely determined
by $(v_j,v_i)=0$ for $i\neq j$ and
$$(v_{j+1},v_{j+1})=\al_{j+1,j}(v_j,v_j).$$
\end{lemma}

$Proof.$ Let us first consider a vector space $V$ with a basis denoted by $(\tilde v_j)$.
We define an action of $B_1$ and $B_2$ on $V$ by substituting $v_j$ by $\tilde v_j$  in the claim, i.e. by $B_1\tilde v_j=[\la-j]\tilde v_j$ and by
$$B_2\tilde v_j= \tilde v_{j+1}+\al_{j-1,j}\tilde v_{j-1}.$$
It is straightforward to check that this action indeed defines a
representation of $U'_q\so_3$; just apply both sides of the given
relation to a basis vector $\tilde v_j$. It also follows directly
that the map $b\mapsto bv_0$ factors over the ideal
$I_\la$ of $\in U'_q\so_3$. Hence we obtain a map from $\Vl$ onto $V$ which maps $v_j$
to $\tilde v_j$. This shows that the $v_j$ are linearly independent.
As $B_2^jv_0=v_j+\sum_{i=0}^{j-2} \beta_iv_i$, it follows that also
the vectors  $B_2^jv_0$ are linearly independent. If $\la$ is a
half-integer, one checks easily that $v_{2\la +1}$ generates an
ideal spanned by the vectors $v_j$ with $j\geq 2\la +1$. As $\Vl$
has a basis of weight vectors, the maximality of this ideal follows
from a well-known standard argument.

To prove the statement about eigenvalues, we use the representations
of $U'_q\so_3$ in \cite{WSp}. They are given by mapping $B_1$ to
$B\otimes1$ and $B_2$ to $1\otimes B$, where $B\in\End(S^{\ot 2})$
and $1$ stands for the identity of $S$, with $S$ the spinor
representation as described in previous sections. It is well-known
that $B_1$ and $B_2$ are conjugated via certain braiding morphisms,
and these braiding morphisms are in the algebra generated by $B_1$
and $B_2$ (see Section \ref{quantugroupeigs}).

 Let $(\ ,\ )$ be an invariant
form on $\Vl$. If $q$ is not a root of unity, then $[\la-j]\neq
[\la-i]$ for $i\neq j$. Hence, by invariance, the $v_j$ are
pairwise orthogonal. But then we also have
$$(v_{j+1},v_{j+1})=(B_2v_j-\al_{j-1,j}v_{j-1},v_{j+1})=(v_j,B_2v_{j+1})=
(v_j, v_{j+2}+\al_{j,j+1}v_j).$$
The claim now follows from the fact that $(v_{j-1},v_{j+1})=(v_j,v_{j+2})=0$.\qed

\medskip
\begin{lemma}\label{sothreeunity} Let $q$ be a primitive $2\ell$-th root of unity, and let $0\leq \la\leq \ell/2$, with $\la$ being a half-integer.
Then there exists a unique simple unitary $U'_q\so_3$ module with highest weight $\la$.
\end{lemma}

$Proof.$ The proof goes along the lines of Lemma \ref{sothreeVerma} by showing that any module as in the statement induces a unique form on $\Vl$.
The main problem now is that $B_1$ has large eigenspaces on $\Vl$.
First assume $\la<\ell/2$. Then we can construct vectors $v_j$, $0\leq j\leq 2\la +1$
with the same inner products as before. In particular,
we have $(v_{2\la+1}, v_{2\la+1})=0$. As the pullback of the form $(\ ,\ )$ on
$\Vl$ is positive semidefinite, it follows that $v_{2\la+1}$ is in its annihilator ideal.
Hence also the vectors $\tilde v_{2\la+1+j}=B_2^jv_{2\la+1}$
are in the annihilator ideal. As the vectors $v_j$ respectively $\tilde v_j$ are of the form $B_2^jv_0 + $ {\it lower terms},
the form is uniquely determined on $\Vl$.

The same strategy also works for $\la=\ell/2$ until the construction
of $v_{\ell}$. We know from the generic case that, in $\Vl$, we have
$v_{2\la+1}=\prod_{j=0}^{2\la} (B_2-[\la-j])v_0$, see Lemma
\ref{sothreeVerma}. As $B_2$ acts via a diagonalizable matrix in a
unitary representation $W$, $v_{2\la +1}$ must be in the annihilator
ideal of the pull-back of the positive definite form on $W$. So, in
particular, also $(v_{\ell+1},v_{\ell-1})=0$ if $\la=\ell/2$. Using
this, we can prove the claim as before for $\la<\ell/2$. \qed

\subsection{$\so_4$ and $\so_5$} First recall the weight structures
for Verma modules for $\so_4$. We have seen in the last subsection
that there exist
polynomials $P_j$ of degree $j$ such that $v_j=P_j(B_2)v_0$ is a
weight vector of weight $\la-j$, where $v_0$ is the highest weight
vector of the Verma module of $U'_q\so_3$ with highest weight $\la$.
Then also $B_3^kP_j(B_2)v_\la$ is an eigenvector of $B_1$ with
eigenvalue $[\la_1-j]$, where $v_\la$ is the highest weight vector
of a $U'_q\so_4$ highest weight module. In view of Lemma \ref{weightshift},
it follows by induction on $j$ that the eigenvalues of $B_3$ are of
the form $[\la_2-j+2i]$, $0\leq i\leq j$. This can be written as

$$\prod_{i=0}^j(B_3-[\la_2-j+2i])P_j(B_2)v_\la=0.$$
Now leaving out the factor for a fixed $i=i_0$ gives us a weight
vector of weight $(\la_1-j,\la_2-j+2i_0)$, or, possibly the zero
vector. As $(\la_2-1,\la_1+1)$ and $(-\la_2-1,-\la_1-1)$ are not
weights of the simple $U'_q\so_4$ module with highest weight
$\la=(\la_1,\la_2)$, the just mentioned expressions for these
vectors have to be in an ideal of the Verma module. This means they
are in the annihilator ideal of any invariant form in the generic
case. Indeed, it follows from Harish-Chandra's theorem (see e.g.
\cite{V}, Theorem 4.7.3) that these vectors generate the maximal
ideal in the classical case. In view of our explicit basis, this can
also be checked directly for $U'_q\so_4$ in the generic case.

If $q$ is a primitive $2\ell$-th root of unity, and $0\leq\la_2\leq
\la_1\leq \ell/2$, it is straightforward to check that the weight
vectors mentioned in the last paragraph are also in the annihilator
ideal of any invariant form, using Lemma \ref{bilinunique}, except
possibly if $\la_1=\ell/2$ and $|\la_2|$ is equal to $\ell/2$ or
$\ell/2-1$. In the first case, we basically have a $U'_q\so_3$
module, as, e.g. for $\la_2=\ell/2$ we have
$B_3B_2v_\la=[\la_2-1]B_2v_\la$ and the claim follows from the
previous section. Similarly, if $\la_2=\ell/2-1$, one considers the
quotient of $\Vl$ modulo the vector of weight $(\ell/2-2,\ell/2+1)$.
It is not hard to check that it is the sum of two $U'_q\so_3$
modules with highest weights $\ell/2$ and $\ell/2-1$, and the claim
again follows from Lemma \ref{sothreeunity}. We have shown most of
the following lemma:

\begin{lemma}\label{sofourunique} Let $q=e^{\pm \pi i/\ell}$.
There is at most one simple unitary $U'_q\so_4$ module with highest
weight $\la$ for any restricted dominant weight $\la$.
The same uniqueness statement holds for a unitary $U'_q\so_5$ module with highest weight $\la=(\ell/2,\ell/2)$ or
$\la=(\ell/2,\ell/2-1)$, provided its restriction to $U'_q\so_4$ is isomorphic to the corresponding restriction for
 the $U'_q\so_5$ module  in Corollary \ref{fusion2} with
the same highest weight $\la$.
\end{lemma}

$Proof.$ After the previous discussion, it only remains to check the claim for
the two $U'_q\so_5$ modules. This can be done by a straightforward inspection as follows: One first checks that all the inner products for $U'_q\so_4$ highest weight vectors are uniquely determined by the value of $(v_\la,v_\la)$,
by Lemmas \ref{bilinunique} and \ref{sothreeunity}. To do this, one deduces from the character formulas in Lemma \ref{comparison} and Theorem \ref{duality} that for $\la=(\ell/2,\ell/2)$, the corresponding $U'_q\so_5$ module decomposes as a direct sum of simple $U'_q\so_4$-modules with highest weights $(\ell/2,j)$ and highest weight vectors $P_j(B_4)v_\la$,
for which the inner products are known by Lemma \ref{sothreeunity}. The same method works for $\la=(\ell/2,\ell/2-1)$,
except for the submodules with highest weights $(\ell/2-1,\pm(\ell/2-1))$. In the latter exceptional cases, the uniqueness of
the norm can be deduced using Lemma \ref{bilinunique}.
The claim now follows from this and  and the already proven claim for unitary $U'_q\so_4$ modules. \qed

\section{Quantum torus and braid representations}\label{section: quantum torus}

\subsection{Quantum torus}
Let $n>1$ and let $A$ be an $(n-1)\times (n-1)$ integer matrix
defined by $a_{ij}=(j-i)$ if $|i-j|=1$ and by $a_{ij}=0$ otherwise.
The quantum $(n -1)$-torus associated with $A$  is:
$$T_q(n) := \C\langle u_1^{\pm 1},\ ...\
u_{n-1}^{\pm 1},\
 : u_iu_j = q^{a_{ij}}u_ju_i\rangle.$$

For $q\in\C^*$ we may specialize $T_q(n)$ at $q$.  In this situation we can
give $T_q(n)$ the structure of a $*$-algebra by setting $u_i^*=u_i^{-1}$.

We have the following elementary lemma:
\begin{lemma}
The algebra $T_q(n)$ has a basis consisting of the monomials
$u_1^{m_1}u_2^{m_2}\ ...\ u_{n-1}^{m_{n-1}}$ with
$m_j\in\Z$ for $1\leq j<n$.
\end{lemma}

$Proof.$ The spanning property is easy to check, using the fact that the generators
$u_i$ commute up to multiplication by a power of $q$. To prove linear independence,
we define an action of $u_i$ on the space of Laurent polynomials $\C[x_1^{\pm 1}, x_2^{\pm 1},\ldots, x_{n-1}^{\pm 1}]$ by
$$u_ix^{\vec m}=q^{m_{i-1}}x_ix^{\vec m},$$
where $\vec m\in \Z^{n-1}$ and $x^{\vec m}=x_1^{m_1}x_2^{m_2}\ ... x_{n-1}^{m_{n-1}}$.
We leave it to the reader to
check that this is indeed a representation of $T_q(n)$. The linear independence follows
from $u^{\vec m}1=x^{\vec m}$ and the linear independence of the vectors $x^{\vec m}$.\qed
\bigskip

\subsection{Finite dimensional representations}\label{fdreps}
If $(\rho,V)$ is a $d$-dimensional representation of $T_q(n)$ for $n\geq 3$ then
$u_iu_{i+1}u_i^{-1}=qu_{i+1}$ implies that $\Spec(\rho(u_i))$ is invariant
under multiplication by $q$.  This, in turn, implies that $q^k=1$ for some $k$
dividing $d$.
Moreover, it is easy to check that $q^k=1$ if and only if $u_i^k$ is in the center
of $T_q(n)$. We define for any $\vec z\in S^{n-1}$, where $S=\{ z\in \C,
|z|=1\}$, the quotient $T^k_q(n,\vec z)$ of $T_q(n)$ (specialized at a primitive $k$th root of unity) via the additional
relations $u_i^k=z_i^k$, $1\leq i\leq n-1$.

\begin{proposition} (a) The algebra $T_q(n)$ has nontrivial finite dimensional
representations if and only if $q$ is a root of unity of finite order.

(b) The algebra $T^k_q(n,\vec z)$ has dimension $k^{n-1}$. It has one simple
module of dimension $k^{(n-1)/2}$ for $n$ odd, and $k$ non-isomorphic simple
modules of dimension $k^{(n-2)/2}$.
\end{proposition}

$Proof.$ Part (a) has been proved already. It also follows easily
that the dimension of $T_q(n,\vec z)$ is at most as stated in (b).
To prove the remainder of (b),
 suppose first that $n$ is odd so that $T^k_q(n,\vec z)$ has an
even number of generators: $u_1,\ldots,u_{n-1}$. Let $V$ be a
$k^{(n-1)/2}$-dimensional vector space with basis $v(\vec i)$, where
$\vec i\in \{0,1,\ldots, k-1\}^{(n-1)/2}$. The action of $u_{2s-1}$
on $V$ is defined by $u_{2s-1}v(\vec i)=z_{2s-1}q^{i_s}v(\vec i)$.
The action of $u_{2s}$ is given  by the rule (indices modulo $k$):

$$u_{2s}(v(i_1,\ldots,i_{s},i_{s+1},\ldots,i_{\frac{n-1}{2}}))=z_{2s}v(i_1,\ldots,
i_s+1 , i_{s+1}-1,i_{s+2},\ldots,i_{\frac{n-1}{2}});$$
in other words, the even indexed generators $u_{2s}$ permute the vectors
$v(i_1,\ldots,i_{\frac{n-1}{2}})$ by shifting the $s$th index up by 1 and the
$(s+1)$th index down by 1, except for $s=(n-1)/2$ where there is no  index
left for shifting down.

It is straightforward to check that  $V$ is a $T^k_q(n,\vec z)$-module.
Standard arguments show that if $W$ is a submodule of $V$, it must contain at least
one common eigenvector of the elements $u_{2s-1}$, $1\leq s<n/2$, i.e. one
of our basis vectors. It then follows for $n$ odd that $W$ contains all basis vectors,
i.e. $W=V$ is simple.  It follows that the image of $T_q(n,\vec z)$ is the full matrix
ring on $V$. This proves all the statements in (b) for $n$ odd.

For $n$ even, we look at the restriction of the just constructed representation
of $T^k_q(n+1,\vec z)$ to $T^k_q(n,\vec z)$. It obviously must be faithful.
On the other hand, it decomposes into the direct sum of $V_r$, $0\leq r<k$
of $T^k_q(n,\vec z)$-modules, where each $V_r$ is the span of
vectors $v(\vec i)$ for which the sum of the  indices $i_1+i_2+\cdots+i_{(n-1)/2}$
is congruent to $r$ mod $k$. From this follow the remaining statements of (b)
for $n$ even. \qed

In what follows we will only need to deal with the special case $\vec z=(1,\ldots,1)$ for which we set
$T^k_q(n)=T^k_q(n,(1,\ldots,1))$.

\subsection{$U'_q\so_n$ representations into the quantum torus} Let $B_i, 1\leq
i<n$ be the generators of $U'_q\so_n$, as before.

\begin{lemma}\label{qtorusrep}
\begin{enumerate}
\item[(a)]  The assignments $B_i\rightarrow \pm\frac{u_i-u_i^{-1}}{q-q^{-1}}$
and
$B_i\rightarrow \pm \i\frac{u_i+u_i^{-1}}{q-q^{-1}}$ extend to algebra
homomorphisms $U^\prime_q\mathfrak{so}_{n}\rightarrow T_{q}(n)$ (for arbitrary
$q$).

\item[(b)] For $q=e^{2\pi \i/(2N)}$ the assignments in (a) extend to algebra
homomorphisms $U^\prime_q\mathfrak{so}_{n}\rightarrow T^{2N}_{q}(n)$.

\item[(c)] (Even case: $N=2k$) Denote by
$\Psi:U^\prime_q\mathfrak{so}_{n}\rightarrow T^{2N}_{q}(n)$
the algebra homomorphism determined by $B_i\rightarrow
b_i:=\i\frac{u_i+u_i^{-1}}{q-q^{-1}}$ as in (b). Then
$\prod_{j=-\frac{N}{2}}^{\frac{N}{2}} (B_i-[j])\in\ker(\psi_N)$
where $[j]:=\frac{q^j-q^{-j}}{q-q^{-1}}$.  Moreover, the set
$[j]$ are distinct for $-\frac{N}{2}\leq j\leq \frac{N}{2}$.

\item[(d)] (Odd case: $N=2k+1$) Denote by
$\Psi:U^\prime_q\mathfrak{so}_{n}\rightarrow T^{2N}_{q}(n)$ the map
determined
by $B_i\rightarrow b_i:=\i\frac{u_i+u_i^{-1}}{q-q^{-1}}$ as in (b).  Then
$\prod_{j=-k-1}^{k} (B_i-[j+\frac{1}{2}])\in\ker(\Psi)$.  In particular the
image of the subalgebra of $U^\prime_q\so_n$ generated by $B_i^2$ factors
through the algebra $Uo_q(n,k)$ (see \cite[Definition 4.7(c)]{WSp}), so that
$\Psi$ induces $\widehat{\psi}_N:Uo_q(n,k)\rightarrow T^{2N}_{q}(n)$.
Moreover, the $b_i$ eigenvalues $[j+\frac{1}{2}]$ for $-k-1\leq j\leq k$ and
the $b_i^2$ eigenvalues $[j+\frac{1}{2}]^2$ for $0\leq j\leq k$  are distinct.

\end{enumerate}
\end{lemma}

$Proof.$
Part (a) is a straight-forward calculation: the case $n=3$ is sufficient since
far-commutation is obvious, and writing out the $q$-Serre relations with
$B_i=x(u_i\pm u_i^{-1})$ gives the specified values of $x$.

Part (b) is obvious since $q^2\neq 1$.

For (c) first observe that for $-\frac{N}{2}\leq j\leq \frac{N}{2}$ the
$N+1$ numbers $[j]$ are distinct since $\sin(x)$ is increasing on $[-\pi,\pi]$.
We have $u_i^{2N}=1$ so $\Spec(u_i)=\{q^j:0\leq j\leq 2N-1\}$, where
$q=e^{\pi \i/N}$ and $\i=q^{\frac{N}{2}}$.  Thus
$\Spec(b_i)=\{q^{N/2}\frac{(q^j+q^{-j})}{q-q^{-1}}:0\leq j\leq
2N-1\}$.  Since $q^{N/2}q^{-j}=-q^{-j-N/2}$ and
$$\{j+N/2\pmod{2N}:0\leq
j\leq 2N-1\}=\{j\pmod{2N}:0\leq
j\leq 2N-1\}$$
we have $\Spec(b_i)=\{[j]:-N/2\leq j\leq N/2\}$.
Since $b_i^*=-b_i$ the minimal polynomial of $b_i$ is a product of distinct
(linear) factors.

For (d) we note as above that $\{[j+1/2]:-k-1\leq j\leq k\}$ is a set of $2k+2$
distinct numbers and $\{[j+1/2]^2:0\leq j\leq k\}$ is a set of $k+1$ distinct
numbers (since $[j+1/2]=-[-j-1/2]$).  We have $u_i^{2N}=1$ so
$\Spec(u_i)=\{q^{j}:0\leq j\leq 2N-1\}$
(where $q=e^{2\pi \i/(2N)}$) and $\i=q^{N/2}$.  Thus
$\i(q^j+q^{-j})=q^{j+N/2}-q^{-j-N/2}$, and $\{j+N/2 \pmod{2N}:0\leq
j\leq 2N-1\}=\{j+1/2\pmod{2N}:0\leq j\leq 2N-1\}$
so we have $\Spec(b_i)=\{[j+1/2]:-k-1\leq j\leq k\}$.
As in (c), the minimal polynomials of $b_i$ and $b_i^2$ are products of
distinct linear
factors. \qed
\bigskip

\subsection{Basics from subfactor theory} In order to compare the representations defined in this section with the ones defined before in connection with fusion categories
we shall need a few basic results from Jones' theory of subfactors (see \cite[Section 3.1]{Jo}).
Let $\A\subset \B$ be finite or infinite dimensional unital von
Neumann algebras with the same identity. Assume that $\B$ has a
finite trace $tr$ satisfying $tr(1)=1$ and $(b,b)=tr(b^*b)>0$ for
$b\neq 0$. Let $L^2(\B,tr)$ be the Hilbert space completion of  $\B$
under the inner product $(\ ,\ )$, and let $e_A$ be the orthogonal
projection onto $L^2(\A,tr)\subset L^2(\B,tr)$. It can be shown that
it maps any element $b\in\B$ to an element $\epsilon_\A(b)\in\A$.
The algebra $\langle \B, e_A\rangle$ is called Jones' basic
construction for $\A\subset \B$. If $\A\subset\B\subset \Ca$ are
finite dimensional algebras and $e\in\Ca$ is such that
$ebe=\epsilon_\A(b)e$ for  all $b\in\B$ and the map $a\mapsto
ae$ defines an algebra isomorphism between $\A$ and $\A e$, one can
show that $\langle \B,e\rangle \cong \B e\B\oplus \B'$, where $\B
e\B$ is isomorphic to a Jones basic construction for $\A\subset\B$,
and $\B'$ is isomorphic to a subalgebra of $\B$.

\subsection{Algebra isomorphisms} We consider the following set-up:
Let $A_i$, $i\in\N$ be a sequence of self-adjoint operators acting on a Hilbert space, satisfying the following conditions:
\begin{enumerate}
\item[(1)]We have $[A_i,A_j]=0$ for $|i-j|>1$, and $\A_{i,j}=\langle A_i, A_{i+1},\ldots, A_{j-1}\rangle$ is a finite-dimensional algebra
for all $i<j$.
\item[(2)]The map $A_i\mapsto A_{i+1}$, $1\leq i\leq j-2$ induces an isomorphism between $\A_{1,j-1}$ and $\A_{2,j}$.
\item[(3)] There exists a unital trace on the algebra $\A$ generated by the elements $A_i$, $i\in\N$,
and an $m>0$ such that $\A_{i,j+1}=\langle A_{i,j},e_{j}\rangle$ is isomorphic to a
Jones basic construction for $\A_{i,j-1}\subset \A_{i,j}$
whenever $j-i\geq m$, where $e_i$ is an eigenprojection of
$A_i$.
\end{enumerate}

\begin{remark}\label{checkconditions} The conditions above are
satisfied for any self-dual object $X$ in a braided unitary fusion category
for which $\End(X^{\otimes 2})$ is generated by an element $A$ for
which $e$ is the projection onto the trivial object $\1\subset
X^{\otimes 2}$, and where $\End(X^{\otimes n})$ is generated by the
elements $A_i=1_{i-1}\otimes A\otimes 1_{n-1-i}$ (see e.g. \cite{WSp}, Prop 2.2 and
the references given there).
\end{remark}

\begin{lemma}\label{algextensions} Let $(A_i)$ and $(\tilde A_i)$ be operators satisfying
conditions $(1)-(3)$ above, with $m\in\N$ as in condition (3). Assume
that $\Phi:A_i\mapsto \tilde A_i$, $1\leq i\leq m$ defines an
algebra isomorphism between $\A_{1,m+1}$ and $\tilde\A_{1,m+1}$.
Then we can extend $\Phi$ to an algebra isomorphism between
$\A_{1,\infty}$ and $\tilde\A_{1,\infty}$ such that $A_i$ is mapped to
$\tilde A_i$ for $i\geq m$. We call this an inclusion-respecting
isomorphism between these algebras.
\end{lemma}

$Proof.$ It follows from our conditions that we can extend $\Phi$ to
an algebra isomorphism between $\A_{1,\infty}$ and
$\tilde\A_{1,\infty}$ by mapping $e_i$ to $\tilde e_i$ for $i>m$, by
uniqueness of the basic construction. It remains to show that it maps
$A_i$ to $\tilde A_i$ for $i>m$. We show this for the algebras $\A_{1,j}$  by
induction on $j$, with $j\leq m+1$ established by assumption. For
the induction step $j\to j+1$, we extend $\Phi$ to $\A_{1,j+1}$ by
mapping $e_j$ to $\tilde e_j$. This also defines an injective
homomorphism from $\A_{2,j+1}$ into the algebra generated by $\tilde
\A_{2,j}$ and $\tilde e_j$, which is a subalgebra of $\tilde
\A_{2,j+1}$. By injectivity and dimension count, the image actually
is $\tilde\A_{2,j+1}$.

On the other hand, using the induction assumption and the isomorphisms
of condition (2), there exists an isomorphism between $\A_{2,j+1}$
and $\tilde\A_{2,j+1}$ which maps $A_i$ to $\tilde A_i$ for $2\leq
i\leq j$. As it also maps $e_j$ to $\tilde e_j$, it must coincide
with the restriction of $\Phi$ to $\A_{2,j}$. This shows the claim. \qed

\subsection{Identifying the representations} We use the notation
$\ve=(1/2,1/2,\ldots, 1/2)\in
\R^j$ and $\e_i$ for the $i$-th standard basis vector of $\R^j$. We associate
these vectors with weights of $\so_n$ for $n=2k$ or $n=2k+1$ in the usual way.

\begin{theorem}\label{centralizercor} We have the following inclusion-respecting isomorphisms
(in the sense of Lemma \ref{algextensions}) where
$\Psi:U'_q\so_n\rightarrow T^{2N}_q(n)$ is as in Lemma
\ref{qtorusrep} and $\Phi:U'_q\so_n\rightarrow \End(S^{\ot n})$ for $N$ even, respectively $\Phi:U'_q\so_n\rightarrow \End(\tilde{S}^{\ot n})$:
\begin{enumerate}
\item[(a)] For $N$ even, $\End(S^{\ot n})=\Phi(U^\prime_q\so_n)\cong\Psi(U'_q\so_n)=\langle 1, u_1+u_1^{-1},\ldots, u_{n-1}+u_{n-1}^{-1}\rangle$.
\item[(b)] For $N$ odd, $\End(\tilde{S}^{\ot n})\supset\Phi(U^\prime_q\so_n)\cong\Psi(U'_q\so_n)=\langle 1, u_1+u_1^{-1},\ldots, u_{n-1}+u_{n-1}^{-1}\rangle$.
\item[(c)] For $N$ odd, $\End(S^{\ot n})\cong\Psi(\langle 1,B_1^2,\ldots,B_{n-1}^2\rangle)=\langle 1, u_1^2+u_1^{-2},\ldots,u_{n-1}^2+u_{n-1}^{-2}\rangle$, where $S\in SO(N)_2$ is the fundamental spinor object.
\end{enumerate}
\end{theorem}

$Proof.$  Parts (a) and (b) are proved by checking that conditions (1)-(3) of
Subsection \ref{fdreps} and Lemma \ref{algextensions} are satisfied for $A_i=\Phi(B_i)$ and for $\tilde A_i=\Psi(B_i)$.
Conditions (1) and (2) are easy to check, using Remark
\ref{checkconditions} and the fact that $u_i\mapsto u_{i+1}$ also
induces a homomorphism in the quantum torus with $q$ a root of
unity.  Indeed, $\tilde S$ is a self-dual object in  $O(N)_2$ and the element $A_1\in\End(S^{\ot 2})$ generates the image of $\Phi(U^\prime_q\so_2)$.

Observe that the representation $\Psi$ of $U'_q\so_n$ into
$T^{2N}_q(n)$ for $q=e^{2\pi \i/(2N)}$ in the previous section has the
same simple components (though not with the same multiplicities) as
its representation $\Phi$ into $\End(S^{\otimes n})$ respectively
$\End(\tilde S^{\otimes n})$ in Corollary \ref{fusion2} for $n\leq 5$.
Indeed, for $n=2$ it suffices to calculate the eigenvalues of $B_1$,
which was done in Lemma \ref{qtorusrep}. They coincide with the ones
in the fusion representation, see \cite[Lemma 4.2 and Proposition 4.3]{WSp}.
It is now easy to check that the usual trace for the
standard representation of the quantum torus satisfies the same
conditions as the functions $\chi^\rho_n$ of Lemma
\ref{rhocharacter}. Hence the same irreducible characters of
$U'_q\so_n$ for $n$ even appear in its representation into the
quantum torus as in its representation into $\End(S^{\otimes n})$
respectively $\End(\tilde S^{\otimes n})$. But as unitary representations
are uniquely determined by their highest weights for $n\leq 5$,
with the additional condition on the restriction for $n=5$ by
Lemma \ref{sofourunique} (observe that all entries $\mu_i$ of our
weights have absolute value $\leq \ell/2$), thus the irreducible
representations of $U'_q\so_n$ in the quantum torus coincide with
the ones in the fusion category, for $n\leq 5$.

Finally, condition (3) of Subsection \ref{fdreps} holds for the algebras
$\A_{i,j}$ with $m=4$ by Remark \ref{checkconditions} and it was verified by
Jones for the algebras $\tilde\A_{i,j}$, see \cite{jonespjm}.
But now the
conditions of Lemma \ref{algextensions} are satisfied for $A_i$ and $\tilde A_i$ with $m=4$, and parts (a) and
(b) follow.

Now suppose $N$ is odd. Denote by $\Da\subset U^\prime_q\so_n$ the
algebra generated by the $(B_i)^2$.  Clearly the
inclusion-respecting isomorphism of (b) restricts to
$\Phi(\Da)\cong\Psi(\Da)$.  Now it is an easy exercise in computing Bratteli diagrams (cf. \cite[Section 5]{jstatmech}) to see that $\dim\Psi(\Da)=\dim\End(S^{\ot n})$.  It follows from this and Remark \ref{nosurjrem}(b) that
$\Phi(\Da)\cong\End(S^{\ot n})$.\qed

\subsection{Braid representations into quantum torus} The
isomorphism in the last theorem transports the braid representations
from the fusion categories to braid representations into the quantum
torus. We determine precisely the images of the braid generators in
these representations, up to an overall scalar factor.

\begin{proposition}\label{identifiedBnreps} Let $q=e^{\pi \i/N}$ and $\psi:\B_n\rightarrow T^{2N}_q(n)$
the braid group representations obtained as compositions of
$\gamma_S:\B_n\rightarrow \Aut(S^{\ot n})$ from Subsection \ref{quantugroupeigs} and the isomorphisms of Theorem \ref{centralizercor}(a,c).
Then:
\begin{enumerate}
\item[(a)] For $N$ odd, we
have $\psi(\sigma_i)=\frac{\gamma}{\sqrt{N}}\sum_{j=0}^{N-1}Q^{j^2}u_i^{2j}$ where
$Q=q^2=e^{2\pi \i/N}$ and $\gamma$ is a scalar of norm 1.
\item[(b)] For $N$ even, we have
$\psi(\sigma_i)=\frac{\gamma}{\sqrt{2N}}\sum_{j=0}^{2N-1}x^{\alpha j^2}u_i^j$
where $x=e^{\pi \i/(2N)}$, $\alpha=1-N(-1)^{N/2}$ and $\gamma$ is a scalar of norm 1.
\end{enumerate}
\end{proposition}
$Proof.$ Clearly $\psi(\sigma_i)$ must be a polynomial in
$b_i=\frac{u_i+u_i^{-1}}{q-q^{-1}}$ for $N$ even and $b_i^2$ for $N$
odd.  Since the isomorphisms of Theorem \ref{centralizercor} respect
inclusions it is enough to prove that
$\psi(\sigma_1)=R_o:=\frac{\gamma}{\sqrt{N}}\sum_{j=0}^{N-1}Q^{j^2}u_1^{2j}\in
T^{2N}_q(n)$ for $N$ odd and
$\psi(\sigma_1)=R_e:=\frac{\gamma}{\sqrt{2N}}\sum_{j=0}^{2N-1}\tau(x)^{j^2}u_1^j\in
T^{2N}_q(n)$ for $N$ even (for some scalars $\gamma$ and some
$\tau\in\Aut_\Q(\Q(x))$).  Comparing the coefficients of $u_1^j$ and
$u_1^{-j}$ one sees that $R_o$ and $R_e$ are indeed polynomials in
$b_1^2$ respectively $b_1$.  Since the number of distinct eigenvalues of
$b_1$ and $b_1^2$ is equal to the dimension of $\End(S^{\ot 2})$
(for $N$ even, respectively odd) it is enough to verify that the
eigenvalues of $R_o$ and $R_e$ coincide with those of $c_{S,S}$ in
Lemmas \ref{Neveneigs} and \ref{Noddeigs} on each $B_1$-eigenspace.
The eigenvalues of $B_1$ are computed in \cite[Lemma 4.2]{WSp}: the
eigenvalue of $B_1$ on the projection onto $V_{[1^{N/2-j}]}$ is
$[j]$ (note that in \cite{WSp} the Young diagram in the subscript
has a typo: $2k$ should be replaced by $k=N/2$ as we have here).
For $N$ odd, we must verify that
$R_ov=\i^{(N/2-s)^2+s}e^{\frac{-s^2\pi \i}{2N}}v$ for any
eigenvector $v$ of $b_1^2$ with eigenvalue $[N/2-s]^2$, for $0\leq
s\leq (N-1)/2$ (up to a scalar independent of $s$) and
 for $N$ even $R_ev=\eta(N/2-s)f(N/2-s)v$ for any eigenvector $v$ of $b_1$ with eigenvalue
$[s]$, for $-N/2\leq s\leq N/2$ where $\eta$ and $f$ are functions defined in Lemma \ref{Neveneigs}
(up to a scalar, independent of $s$, and some choice of $\tau$).

We will give the details in the $N$ even case and leave the $N$ odd case to the reader.

For $N$ even and $-N/2\leq s\leq N/2$, $u_1$ acts on the $[s]$-eigenspace of $b_1$ by $x^{\pm(2s-N)}$.
The corresponding eigenvalue of $\frac{1}{\sqrt{2N}}\sum_{j=0}^{2N-1}(x)^{j^2}u_1^j$ is:
$$\frac{1}{\sqrt{2N}}\sum_{j=0}^{2N-1}x^{j^2\pm(2s-N)j}.$$ Completing the square we have:
$$\frac{x^{-(s-N/2)^2}}{\sqrt{2N}}\sum_{j=0}^{2N-1}x^{(j\pm (s-N/2))^2}.$$ Since $x$ is a $4N$th root
of unity and the set of residues modulo $4N$ of $(j\pm (s-N/2))^2$ is the same for $0\leq j\leq 2N-1$ and $2N\leq j\leq 4N-1$
we double the sum to obtain:
$$\frac{x^{-(s-N/2)^2}}{2\sqrt{2N}}\sum_{j=0}^{4N-1}x^{j^2}=\frac{x^{-(s-N/2)^2}(1+\i)}{\sqrt{2}}.$$ using Dirichlet's improvement on Gauss' result
(see e.g. \cite{Davenport}).

Rescaling (independent of $s$) we obtain the eigenvalue
$f(N/2-s)(-\i)^{(N/2-s)}$ for $R_e$ on these spaces. The result now follows by verifying,
for $\alpha=1-N(-1)^{N/2}$, that
$$\frac{[f(N/2-s)(-\i)^{(N/2-s)]^\alpha}}{(\eta(N/2-s)f(N/2-s))}$$ is
independent of $s$.\qed

\begin{remark}
The Gaussian representations of $\B_n$ in $T_n^{2N}(n)$ described in Proposition \ref{identifiedBnreps}(a) go back at least to \cite{GJ} in the case $N$ is odd and were certainly known to Jones in the case $N=3$ in the early 1980s.
In the case $N$ is even these representations seemed to explicitly appear only recently \cite{GR},
in which results of \cite{jonespjm}
are employed, and their properties are studied in some detail.
\end{remark}

As a consequence we can prove (a generalized version of) \cite[Conjecture 5.4]{RW}:
\begin{theorem}
The images of the braid group representations on $\End_{SO(N)_2}(S^{\ot n})$ for $N$ odd and $\End_{SO(N)_2}(S_\pm^{\ot n})$
for $N$ even are isomorphic to images of braid groups in Gaussian representations; in particular, they are finite groups.
\end{theorem}
$Proof.$ In \cite{GR} the Gaussian representations are shown to have finite image. Hence for $N$ odd, the
claim is immediate from Proposition \ref{identifiedBnreps}. For $N$ even the
same analysis implies that the braid group representation on
$\End_{O(N)_2}(S^{\ot n})$ for $N$ even is a finite group.  Since
the forgetful functor $F:O(N)_2=(SO(N)_2)^{\Z_2}\rightarrow SO(N)_2$
is a braided tensor functor and the braiding is functorial we
conclude that the image of the braid group acting on
$\End_{SO(N)_2}(S_\pm^{\ot n})$ is a (finite) subquotient of the
image of the braid group acting on $\End_{O(N)_2}(S^{\ot n})$. \qed

\end{document}